\documentclass[10pt]{amsart}
\allowdisplaybreaks[4]
\usepackage{color}

\definecolor{c20}{rgb}{0.,0.7,0.}
\definecolor{c30}{rgb}{0.,0.,1.}
\definecolor{c40}{rgb}{1,0.1,0.7}
\definecolor{c50}{rgb}{1,0,0}
\definecolor{c60}{rgb}{1,0.9,0.1}

\def\pE#1{\textcolor{c20}{#1}}
\def\pE#1{#1}

\def\cL#1{\textcolor{c50}{#1}}
\def\cL#1{#1}

\def\tbb#1{\textcolor{blue}{#1}}
\def\tbb#1{#1}

\newcommand{\kb}[1]{\boldsymbol{#1}}
\newcommand{\vk}[1]{\kb{#1}}

\newcommand{\ve}{\varepsilon}

\newcommand{\abs}[1]{\left\lvert #1 \right\rvert}

\newcommand{\E}[1]{\mathbb{E}\left\{#1\right\}}

\newcommand{\pk}[1]{\mathbb{P} \left\{ #1 \right\} }

\newcommand{\R}{\mathbb{R}}

\newcommand{\N}{\mathbb{N}}
\newcommand{\inr}{\in \R}
\newcommand{\inn}{\in \N}

\newcommand{\limit}[1]{\lim_{#1 \to   \infty}}

\newcommand{\BQN}{\begin{eqnarray}}
\newcommand{\EQN}{\end{eqnarray}}
\newcommand{\BQNY}{\begin{eqnarray*}}
\newcommand{\EQNY}{\end{eqnarray*}}

\newcommand{\BS}{\begin{sat}}
\newcommand{\ES}{\end{sat}}
\newcommand{\BT}{\begin{theo}}
\newcommand{\ET}{\end{theo}}
\newcommand{\BK}{\begin{korr}}
\newcommand{\EK}{\end{korr}}

\newcommand{\BD}{\begin{de}}
\newcommand{\ED}{\end{de}}

\newcommand{\BIT}{\begin{itemize}}
\newcommand{\EIT}{\end{itemize}}
\newcommand{\BDI}{\begin{description}}
\newcommand{\EDI}{\end{description}}

\newcommand{\BRM}{\begin{remarks}}
\newcommand{\ERM}{\end{remarks}}

\newcommand{\BEL}{\begin{lem}}
\newcommand{\EEL}{\end{lem}}

\newtheorem{theo}{Theorem}[section]
\newtheorem{sat}[theo]{Proposition}
\newtheorem{de}[theo]{Definition}
\newtheorem{lem}[theo]{Lemma}

\newtheorem{korr}[theo]{Corollary}
\newtheorem{remark}[theo]{Remark}
\newtheorem{remarks}[theo]{Remarks}
\newtheorem{prop}[theo]{Proposition}

\newcommand{\nelem}[1]{{Lemma \ref{#1}}}

\newcommand{\netheo}[1]{{Theorem \ref{#1}}}

\newcommand{\prooftheo}[1]{ \textsc{\bf Proof of Theorem} \ref{#1}:}

\newcommand{\prooflem}[1]{\textsc{\bf Proof of Lemma} \ref{#1}:}

\newcommand{\COM}[1]{}

\newcommand{\QED}{\hfill $\Box$}

\def\rw{\rightarrow}

\def\IF{\infty}

\def\Cov{\mathrm{Cov}}

\date{}

\def\oo{(1+o(1))}

\def\LT{\left}
\def\RT{\right}
\def\H{\mathcal{H}}

\def\ooo{(1+o(1))}

\def\rw{\rightarrow}

\def\Piter{\mathcal{P}}

\def\Del{\triangle}

\def\vn{\varepsilon}
\def\Var{\text{Var}}

\def\NN{\mathcal{N}}

\def\FF{\widetilde{\mathcal{H}}}
\def\PP{\widetilde{\mathcal{P}}}
\def\H{\mathcal{H}}

\def\SigS{\widetilde{\sigma} }

\def\upto{\uparrow}
\def\downto{\downarrow}

\begin{document}

\title[Parisian ruin over a finite-time horizon] {On Parisian ruin over a finite-time horizon}

\author{Krzysztof D\c{e}bicki}
\address{Krzysztof D\c{e}bicki, Mathematical Institute, University of Wroc\l aw, pl. Grunwaldzki 2/4, 50-384 Wroc\l aw, Poland}
\email{Krzysztof.Debicki@math.uni.wroc.pl}

\author{Enkelejd  Hashorva}
\address{Enkelejd Hashorva, University of Lausanne\\
B\^{a}timent Extranef, UNIL-Dorigny, 1015 Lausanne, Switzerland
}
\email{Enkelejd.Hashorva@unil.ch}

\author{Lanpeng Ji }
\address{Lanpeng Ji, University of Lausanne\\
B\^{a}timent Extranef, UNIL-Dorigny, 1015 Lausanne, Switzerland
}
\email{Lanpeng.Ji@unil.ch}

\bigskip

\date{\today}
 \maketitle

{\bf Abstract:} For a risk process $R_u(t)=u+ct-X(t), t\ge 0$, where $u\ge 0$ is the initial capital, $c>0$ is the premium rate and
$X(t),t\ge 0$ is an aggregate claim process, we  investigate the probability of the Parisian ruin
\tbb{
\[
\mathcal{P}_S(u,T_u)=\pk{\inf_{t\in[0,S]} \sup_{s\in[t,t+T_u]} R_u(s)<0},
\] with} a given positive constant $S$ and a positive  measurable function $T_u$.
We derive
asymptotic expansion of $\mathcal{P}_S(u,T_u)$, as $u\to\infty$,  for the aggregate claim process
$X$  modeled by  Gaussian processes.
As a by-product, we derive the exact tail asymptotics of the infimum of a standard Brownian
motion with drift over a finite-time interval.



\section{Introduction}
Consider a random process  $\{X(t),t\ge0\}$ which models the aggregate claim process   of an insurance company, i.e., $X(t)$ represents the total amount of claims paid up to time $t$.
\tbb{In a theoretical insurance model
the main object of interest is the so-called surplus process $R_u$, defined by 
\BQN\label{eq:R}
R_u(t)=u+ct-X(t), \ \ \ t\ge0,
\EQN
where $c>0$ models the premium income rate
and
$u\ge 0$ is the initial capital;} see e.g., \cite{EKM97}.
For any $S\in (0,\IF]$, define  the (classical) {\it probability of ruin} during the time period $[0,S]$ as
\BQN
P_S(u):=\pk{\inf_{t\in[0,S]} R_u(t)<0}.
\EQN
We refer to
 \cite{EKM97, HP99, HP08, dieker2005extremes} and references therein for the literature
on properties of $P_S(u)$.

The contributions \cite{DW08, czarna2011ruin}
introduced and studied the Parisian ruin which allows the surplus process to spend a pre-specified time under level zero before ruin is recognized. Initially, Parisian stopping times have been studied by \cite{Chesney} in the context of barrier options in mathematical finance. \\
Let  $T_u$,  depending eventually on the initial capital $u$,
model the pre-specified time which is a positive deterministic function of $u$.
In our setup, the
 {\it probability of Parisian ruin } over the time period $[0,S]$ is defined as
\BQNY
\mathcal{P}_S(u,T_u)=\pk{\inf_{t\in[0,S]} \sup_{s\in[t,t+T_u]} R_u(s)<0}.
\EQNY
Calculation of the probability of Parisian ruin $\Piter_S(u,T_u)$
is  more complex than the calculation of 
$P_S(u)$. When $S=\IF$ and $X$ is modelled by a specific class of L\'{e}vy processes,
 exact formulas for $\Piter_\IF(u,T)$, with $T\in(0,\IF)$ are derived in
\cite{DW08, czarna2011ruin, LCP13}. See also \cite{Palmowski14,Palm14a,Palm14b,Irmina} for some recent developments.

In this paper, we shall \pE{investigate the probability} of  Parisian ruin when the initial
capital becomes large (tends to infinity) and $X$ is modeled by a Gaussian \cL{process}. 
It appears that the qualitative type of the obtained asymptotics
is different \cL{from that} of  the corresponding L\'evy model.
Specifically, if $X$ is a L\'{e}vy  process such that $X(S)$ has a long-tailed distribution, which in view of
\cite{Foss13} means that
there exists some function $h(u), u\ge 0$ satisfying
$\lim_{u\to \IF} \frac{u}{h(u)}= \lim_{u\to \IF} {h(u)}= \IF$
 such that
 \BQN\label{F13}
\pk{X(S)> u+ h(u)} =  \pk{X(S)>u}\ooo, \quad u\rightarrow\IF,
\EQN
then the following proposition holds.
\begin{prop}\label{prop.levy}
Let $S>0$, and $T_u, u\ge 0$ be a positive bounded measurable function.
\pE{If} $X$ is a L\'evy process such that $X(S)$ has a long-tailed distribution, then
\BQN\label{LevyL}
\mathcal{P}_S(u,T_u)  =  \pk{ X(S) > u}\ooo, \quad u\to \IF.
\EQN
\end{prop}
We give the proof of the above proposition in Section 4.
A straightforward application of Proposition \ref{prop.levy}
for
$X$ being an $\alpha$-stable L\'evy process with $\alpha\in(1,2)$,
(i.e., $X(t)\overset{d}=\mathcal{S}_\alpha(t^{1/\alpha}, \beta,0), t>0$, where
$\mathcal{S}_\alpha(\sigma, \beta, d)$ denotes a stable random variable with index of stability $\alpha$,
scale parameter $\sigma$, skewness parameter $\beta$ and drift parameter $d$; see  e.g., \cite{SamT94}),
implies that
\BQNY
\mathcal{P}_S(u,T_u)   =  
 \frac{(1-\alpha)}{\Gamma(2-\alpha)\cos(\pi \alpha/2)}\left(\frac{1+\beta}{2}\right)
 Su^{-\alpha}\ooo,\ \ u\to\IF,
\EQNY
where  $\Gamma(\cdot)$ denotes the Euler Gamma function.

The above restriction that $X(S)$ is long-tailed excludes the classical case that
$X$ is a standard Brownian motion.
Given the importance of the Brownian motion risk process
(see e.g., \cite{Mich98, HP99, DHJ13b}) in this contribution we shall investigate the
asymptotics of the probability of Parisian ruin $\Piter_S(u,T_u)$ with $S\in (0,\IF)$
for large classes of Gaussian risk processes.
It turns out that in opposite to Proposition \ref{prop.levy}, for this model the asymptotics is highly sensitive to
$T_u$. Details are presented in Section 3.

As shown for instance in \cite{HJP13, HJ13, DHJ13b}, the calculation of the  probability of ruin 
over an infinite-time horizon for Gaussian risk processes \pE{raises} interesting theoretical questions for the asymptotic theory of Gaussian processes and related random fields.
Similarly, the calculation of the probability of Parisian ruin over finite-time horizon (which is more involved than the investigation of the \cL{infinite-time} horizon) \pE{raises several} interesting questions as well. For instance, for our investigations it is crucial to obtain certain extensions of Piterbarg lemma, which we shall present in \nelem{LemGP} in Appendix.
For details on Piterbarg and  Pickands lemmas see e.g., \cite{debicki2002ruin,SheppKL,DebKo2013,HJ13b}.
Another interesting problem motivated by this paper is the investigation of the asymptotic
behaviour of
\[
\pk{\inf_{t\in[T_1,T_2]}\Bigl(X(t)-ct\Bigr)>u}, \quad T_2>T_1>0
\]
as $u\to\infty$ with  $X$  a centered  non-stationary Gaussian process. This problem seems to be very hard; we are able to derive an explicit result only for $X$
being a standard Brownian motion $\{B_1(t),t\ge 0\}$; see Theorem \ref{l.bm}.

This paper is organized as follows: After some preliminary results
\tbb{given in the next section, in Section 3} we present our main findings. \netheo{ThmGauExact} \tbb{provides}
the exact asymptotics of $\mathcal{P}_S(u,T_u)$ for $T_u$ converging to 0.
When $X$ is a standard Brownian motion our result holds for $u^2 T_u \to T\in [0,\IF)$ as $u\to \IF$.
The case of constant or general bounded $T_u$ is investigated in \netheo{ThmGauBou},
which gives an asymptotic lower bound for $\mathcal{P}_S(u,\cL{T})$,
and \cL{in} Theorem \ref{th.log} which provides  logarithmic asymptotics of $\mathcal{P}_S(u,T_u)$.
The proofs of the main results are displayed in Section 4, followed  by an Appendix (Section 5).

\section{Preliminaries}
\def\sigmaX{\sigma}

 Let $\{X(t),t\ge0\}$ be a centered Gaussian process with almost surely (a.s.) continuous sample paths and variance function $\sigmaX^2(\cdot)$.
In our setup, $\sigmaX( \cdot)$ is not a constant function,  and therefore the stationary Gaussian processes are excluded. The theory of extremes of non-stationary Gaussian processes is established in numerous contributions; see e.g.,
\cite{Fatalov92, Pit96}.
\tbb{A} key condition in the case of processes with non-constant variance is its local \tbb{structure} at the maximum point of the variance function;
for our setup we shall assume \tbb{the following local} condition:

 {\bf Assumption A1}.  The standard deviation function $\sigmaX(\cdot)$ of the Gaussian process $ X $ attains
its maximum $\SigS$ on  $[0,S]$ at the unique point $t=S$.  Further, there exist   positive constants
$ \beta_1, \beta_2, A$, and  $A_+>0$ (or $A_-<0$) such that
\BQNY\label{eq:sigmaX1}
\sigmaX(t) =  \SigS -A(S-t)^{\beta_1}\ooo, \quad t\upto S
\EQNY
and
\BQN\label{eq:sigmaX2}
\sigmaX(t) =  \SigS -A_{\pm}(t-S)^{\beta_2}\ooo, \quad t\downto S.
\EQN
It is worth noting that in our setup the behaviour of $\sigmaX(\cdot)$ in the right neighborhood of $S$ can be different from that in the left-neighbourhood of $S$.
Specifically, in condition \eqref{eq:sigmaX2} the constant $A_{\pm}$ can be positive or negative, and moreover
the index $\beta_2$ can be different from the index $\beta_1$.

Our next two assumptions are standard, see Chapter 1 in \cite{Pit96}.

 {\bf Assumption A2}.  There exist  some positive constants
$\alpha\in (0,2], D$ such that
\BQNY\label{eq:covX}
 \Cov\LT(\frac{X(t)}{\sigmaX(t)}, \frac{X(s)}{\sigmaX(s)}\RT)=1- D|t-s|^{\alpha}\ooo, \quad  t,s \rightarrow S. 
 \EQNY

{\bf Assumption A3}. There exist some positive constants $ Q$, $\gamma$ and $S_1<S$ such that, for all $s,t\in[S_1, S]$
\BQN\label{eq:regu}
\E{(X(t)-X(s))^{2}} \leq  Q|t-s|^{\gamma}.
\EQN

Next, we introduce some generalizations of  the Pickands and  Piterbarg constants.
We refer to \cite{PicandsA, Pit96, TY15} for the \cL{definitions} and properties of the
(classical) Pickands  and Piterbarg constants.  See also \cite{DikerY} for alternative formulas of Pickands constant. \\
Let $\{B_\alpha(t),t\inr\}$ be a standard fractional Brownian motion (fBm)
with Hurst index $\alpha/2\in(0,1]$, i.e.,  it is a centered Gaussian process with a.s. continuous sample paths and covariance function
$$
\Cov(B_\alpha(t),B_\alpha(s))=\frac{1}{2}(\abs{t}^\alpha+\abs{s}^\alpha-\abs{t-s}^\alpha),\ \ s,t\inr.
$$
Define the generalized Pickands constant as
\BQN\label{eq:GH}
\FF_{\alpha}(T)=\lim_{\lambda\to\IF}\frac{1}{\lambda}\FF_\alpha(\lambda,T),\ \ T\ge0, , \alpha\in(0,2],
\EQN
where
\BQNY
\FF_\alpha(\lambda,T)=\E{\exp\LT(\sup_{t\in[0,\lambda]}\inf_{ s\in[0,T]}\LT(\sqrt{2}B_{\alpha}(t-s)-\abs{t-s}^{\alpha}\RT)\RT)}\in(0,\IF),\ \ \lambda,T\ge0,  \alpha\in(0,2].
\EQNY
Further, we define the generalized Piterbarg constant as
\BQN\label{eq:GP}
\PP^{b_1,b_2}_{\alpha,\beta}(T)=\lim_{\lambda\to\IF} \PP^{b_1,b_2}_{\alpha,\beta}(\lambda,T),\ \ T\ge0,\ b_1>0, b_2\inr, \alpha\in(0,2], \beta\ge\alpha,
\EQN
where, for any positive constants $\lambda, \beta, b_1, T\ge0, \alpha\in(0,2]$ and $b_2\inr$
\BQNY
\PP^{b_1,b_2}_{\alpha,\beta}(\lambda,T)=\E{\exp\LT(\sup_{t\in[0,\lambda]}\inf_{ s\in[0,T]}\LT(\sqrt{2}B_{\alpha}(t-s)-\abs{t-s}^{\alpha}-b_1\abs{t-s}^{\alpha}I_{(t>s)}-b_2\abs{t-s}^{\alpha}I_{(t\le s,\alpha=\beta)}\RT)\RT)},
\EQNY
with $I_{(\cdot)}$ the indicator function. Note that both $\FF_\alpha(\lambda,T)$ and $\PP^{b_1,b_2}_{\alpha,\beta}(\lambda,T)$ are well defined since
$$
\E{\exp\LT(\sup_{t\in[0,\lambda]}  \sqrt{2}B_{\alpha}(t)\RT)}<\IF,\ \ \forall \lambda\ge0,
$$
which follows directly from Piterbarg inequality (see Theorem 8.1 in \cite{Pit96}).
As it will be seen from the proof of \netheo{ThmGauExact} below,   both $\FF_{\alpha}(T)$ and $\PP^{b_1,b_2}_{\alpha,\beta}(T)$ defined above are positive and finite.
Note further that the  classical  Pickands constant $\H_\alpha$ equals $\FF_\alpha(0)$ and the classical Piterbarg constant $\H_\alpha^{b_1}$ equals $\PP^{b_1,b_2}_{\alpha,\beta}(0)$.

Finally, we present a  theorem on the asymptotics of the infimum of Brownian motion with
linear drift over a finite-time interval, which will be used in the next section and is of some independent interest.
Hereafter  $\Psi(\cdot)$ denotes the tail distribution function of an $N(0,1)$ random variable
and $\varphi(\cdot)$ \cL{is} its density.

\BT\label{l.bm}
For any $c>0$  and two constants $T_2>T_1>0$  we have
\BQN\label{eqLemBm}
\pk{\inf_{t\in[T_1,T_2]}\Bigl(B_1(t)-ct\Bigr)>u}=
K_{c,T_2-T_1} \frac{T_1}{u} \Psi\left(\frac{u+c T_1}{\sqrt{T_1}}\right)\ooo, \quad u\to\IF,
\EQN
where
\BQN\label{eq:KcT}
K_{c,y}=2\varphi \LT( c \sqrt{y} \RT) \frac{1}{\sqrt{y}} -
 2 c  \Psi\LT(  c \sqrt{y} \RT)>0,\ \ \ y>0.
 \EQN
\ET

\section{Main Results}\label{s.main}
In this section, we present our main results on the asymptotic behaviour of $\Piter_S(u,T_u)$ as $u\to\IF$. It turns out that when $T_u$ does not vanish to 0 as $u\to\IF$, the exact asymptotics is very hard to derive.
For such cases we shall give a lower asymptotic bound and then the logarithmic asymptotics of $\Piter_S(u,T_u)$ for $X$ being with stationary increments.
Finally, in \netheo{ThmGauExact} we show the exact asymptotics of  $\Piter_S(u,T_u)$,
under certain restrictions on the speed of convergence of $T_u$ to 0, for $X$ satisfying {\bf A1--A3}.

\subsection{Logarithmic asymptotics}
\label{s.log}

The following theorem displays  an asymptotic lower bound
of $\Piter_S(u,T)$, which is logarithmically exact for all large $u$. We write below  $V'(t)$ for the derivative of the variance function $\sigmaX^2(t)$ if it exists.
\BT\label{ThmGauBou}
Let $\{X(t),t\ge0\}$ be a centered Gaussian
process with a.s. continuous sample paths,  $X(0)=0$ and stationary increments. If further the variance function $\sigma^2 (\cdot)$ is differentiable, strictly increasing and convex, then for any positive constants $S,T$
\BQN
\Piter_S(u,T)\ge C_{c,\Delta}\frac{\sigmaX^2(S) }{u} \Psi\LT(\frac{ u+c S}{\sigmaX(S)}\RT)\ooo,\ \ u\to\IF,
\EQN
where
$$
C_{c,\Delta}=2\varphi \LT( \frac{c \sqrt{\Delta}}{V'(S)}\RT) \frac{1}{\sqrt{\Delta}} -
 2 \frac{c}{V'(S)}   \Psi\LT( \frac{c \sqrt{\Delta}}{V'(S)}\RT),\ \ \Delta=\sigma^2 (S+T)-\sigmaX^2(S).
$$
\ET
The proof of Theorem \ref{ThmGauBou} is given in Section \ref{p.low.bound}.

The next result constitutes an LDP counterpart of Proposition \ref{prop.levy}.
\BT\label{th.log}
Let $\{X(t),t\ge0\}$ be a centered Gaussian
process with a.s. continuous sample paths,  $X(0)=0$ and stationary increments.
If further the variance function $\sigmaX^2(\cdot)$
is differentiable, strictly increasing and convex, then \cL{for any bounded measurable function $T_u>0$ and any $S>0$}
\BQN\label{log;kr}
\lim_{u\to\infty}\frac{\log(\Piter_S(u,\cL{T_u}))}{u^2}=-\frac{1}{\sigmaX^2(S)}.
\EQN
\ET
Section \ref{p.log} displays the proof of Theorem \ref{th.log}.

We note that, the claim in \eqref{log;kr} matches the logarithmic asymptotics of the classical
ruin probability, i.e.
\[
\lim_{u\to\infty}\frac{\log(\Piter_S(u,\cL{T_u}))}{\log(\Piter_S(u))}=
\lim_{u\to\infty}\frac{\log(\Piter_S(u,\cL{T_u}))}{\log(\pk{X(S)>u})}=
1
\]
and does not depend on the value of the parameter $c$.

\subsection{Exact asymptotics}
The problem of finding the exact asymptotics of
$\Piter_S(u,T_u)$ needs much more precise analysis.
Next, we discuss  the case that $T_u$ is sufficiently small, tending to 0 as $u\to \IF$.

\BT\label{ThmGauExact}
Let $\{X(t),t\ge0\}$ be a centered Gaussian
process satisfying assumptions \textbf{A1-A3} with the parameters therein, and let $T_u$ be a positive \cL{measurable} function of $u$. Assume that $\beta_1\le\beta_2\le 1$.  For any positive constant $S$, we have,  as $u\to \IF$:\\
(i) If $\alpha< \beta_1$ and $\limit{u} T_uu^{2/\alpha}= T\in[0,\IF)$,  then
\BQN
\Piter_S(u,T_u)=\FF_\alpha(D^{\frac{1}{\alpha}} \SigS^{-\frac{2}{\alpha}} T)  \Gamma\LT(\frac{1}{{\beta_1}}+1\RT)D^{\frac{1}{\alpha}} A^{-\frac{1}{{\beta_1}}}   \SigS^{\frac{3}{{\beta_1}}-\frac{2}{\alpha}} u^{\frac{2}{\alpha}-\frac{2}{{\beta_1}}}\Psi\LT(\frac{u+cS}{\SigS} \RT) \ooo.
\EQN
(ii) If $\alpha= {\beta_1}$ and $\limit{u} T_uu^{2/\alpha}= T\in[0,\IF)$,  then 
\BQN
\Piter_S(u,T_u)= \PP_{\alpha,\beta_2}^{A/(D\SigS), A_\pm/(D\SigS)}(D^{\frac{1}{\alpha}} \SigS^{-\frac{2}{\alpha}} T)\Psi\LT(\frac{u+cS}{\SigS} \RT)  \ooo.
\EQN
(iii) If $\alpha>{\beta_1}$ and $\limit{u} T_uu^{2/{\alpha}}=\limit{u} T_uu^{2/{\beta_2}}=0$,  then 
\BQN
\Piter_S(u,T_u)= \Psi\LT(\frac{u+cS}{\SigS} \RT)  \ooo.  
\EQN
\ET
\begin{remark}  Clearly, if $T_u= 0, u\ge0$,
then  $\Piter_S(u,0)$ becomes the classical probability of ruin $P_S(u)$.
Since, as mentioned  above,
$\FF_\alpha(0)=\H_\alpha$ and  $\PP_{\alpha,\beta_2}^{A/(D\SigS), A_\pm/(D\SigS)}(0)=\H_{\alpha}^{A/(D\SigS)}$, the asymptotics
of $P_S(u)$ is retrieved and agrees with \tbb{findings of \cite{Pit96}}.
\end{remark}

Specialized to the case of the fBm risk process, 
\tbb{the above theorem entails the following result.}
\BK\label{CorfBm}
Let $\{X(t),t\ge0\}$ be a standard fBm with Hurst index $\alpha/2\in(0,1]$. For any positive constant $S$,
we have, as $u\to \IF$:\\
(i) If $\alpha\in(0,1)$ and $\limit{u} T_uu^{2/\alpha}= T\in[0,\IF)$,  then 
\BQNY
\Piter_S(u,T_u)=\FF_\alpha(2^{-\frac{1}{\alpha}} S^{-2} T) \alpha^{-1} 2^{1-\frac{1}{\alpha}} S^{\alpha-1}  u^{\frac{2}{\alpha}-2}\Psi\LT(\frac{u+cS}{ S^{\alpha/2}} \RT) \ooo.
\EQNY
(ii) If $\alpha=1$ and $\limit{u}  T_u u^2= T\in[0,\IF)$,  then
\BQNY
\Piter_S(u,T_u)= \PP_{1,1}^{\cL{1,-1}}( 2^{-1} S^{-2}  T)\Psi\LT(\frac{u+cS}{S^{1/2} } \RT)  \ooo.
\EQNY
(iii) If $\alpha\in(1,2]$ and $\limit{u} T_uu^{2} =0$,  then
\BQNY
\Piter_S(u,T_u)= \Psi\LT(\frac{u+cS}{S^{\alpha/2} } \RT)  \ooo.  
\EQNY
\EK

\begin{remark}
The case that $T_u=T>0$ for all $u$ large is much more difficult to deal with
and most probably needs to develop new techniques that allow derivation of the asymptotics of
tail distribution of infimum of a Gaussian process.
\end{remark}

\begin{remark}
\cL{As in \cite{czarna2011ruin, LCP13, DHJ13b} we define
the {\it Parisian ruin time} of  the  risk process $R_u$ by
\BQNY\label{eq:tau2}
\tau_u=\inf\{t\ge T_u: t-\kappa_{t,u}\ge T_u\},\ \ \ \text{with} \ \kappa_{t,u}=\sup\{s\in[0,t]: R_u(s)\ge0\}.
\EQNY
Under the assumptions of Corollary \ref{CorfBm}  it follows along the lines of the arguments in \cite{DHJ13a}
 that 
\BQN
\lim_{u\to\IF}\pk{u^2(S-\tau_u)\le x\big\lvert \tau_u<S}=1-\exp\LT(-\frac{\alpha}{2}S^{-\alpha-1} x\RT)
\EQN
holds for any $x$ positive.}

\end{remark}


\COM{
We conclude this section with a result concerning the asymptotic distribution of the conditional scaled Parisian ruin time. As in \cite{czarna2011ruin, LCP13, DHJ13b},
the {\it Parisian ruin time} of  the  risk process $R_u$ is defined as
\BQN\label{eq:tau2}
\tau_u=\inf\{t\ge T_u: t-\kappa_{t,u}\ge T_u\},\ \ \ \text{with} \ \kappa_{t,u}=\sup\{s\in[0,t]: R_u(s)\ge0\}.
\EQN
Here we make the convention that $\inf\{\emptyset\}=\IF$ and $\sup\{\emptyset\}=0$. Denote by $\sigmaX'(t)$ the left derivative of $\sigmaX(t)$.
\BT\label{ThmRuintime}
Under the assumptions in \netheo{ThmGauExact} with the conditions in cases (i)--(iii),  we have, for any positive constant $S$ satisfying $\sigmaX'(S)\neq 0$ and any $x\ge0$
\BQNY
\lim_{u\to\IF}\pk{u^2(S-\tau_u)\le x\big\lvert \tau_u<S}=1-\exp\LT(-\frac{\sigmaX'(S)}{\widetilde{\sigma}^3} x\RT),\ \ \ u\to\IF.
\EQNY
\ET
}

\section{Proofs}\label{p.levy}
\subsection{Proof of Proposition \ref{prop.levy}}
First, for any $u$ positive
\BQNY
\mathcal{P}_S(u,T_u)&=& \pk{\sup_{t\in[0,S]}\inf_{s\in[t,t+T_u]}\Bigl(X(s)-cs\Bigr)>u}\\
&\le&\pk{\sup_{t\in[0,S]}X(t)>u}.
\EQNY
Further, in view of \cite{Berman86}  we have
 \BQNY
\pk{\sup_{t\in[0,S]}X(t)>u}= \pk{ X(S)>u}\ooo,\ \ u\to\IF
\EQNY
implying thus
\BQNY
\mathcal{P}_S(u,T_u)\le \pk{ X(S)> u}\ooo,\ \ u\to\IF.
\EQNY
We derive next the lower bound.
Taking  $h(\cdot)$ to be such that \eqref{F13} holds we have
\BQNY
\lefteqn{\pk{\sup_{t\in[0,S]}\inf_{s\in[t,t+T_u]}\Bigl(X(s)-cs\Bigr)>u}}\\
&\ge & \pk{\inf_{s\in[S,S+T_u]}\Bigl(X(s)-cs\Bigr)>u}\\
&\ge&\pk{\inf_{t\in[S,S+T_u]}  \Bigl( X(t) - X(S) -c(t-S) + X(S) -cS\Bigr) >u, X(S) -cS>u+h(u)}.
\EQNY
 Since $T_u$ is bounded, we have $\sup_{u\in[0,\IF)}T_u<M$ for some constant $M$. By the fact that $X$ has independent and stationary increments we may further write
\BQNY
\lefteqn{\pk{\sup_{t\in[0,S]}\inf_{s\in[t,t+T_u]}\Bigl(X(s)-cs\Bigr)>u}}\\
&\ge&\pk{\inf_{t\in[0,M]} \Bigl(X(t)-ct \Bigr)>-h(u)}\pk{ X(S) -cS >u+h(u)}\\
&=& \pk{X(S)> u} \ooo,\ \ u\to\IF
\EQNY
establishing the proof. \QED

\subsection{Proof of Theorem \ref{l.bm}}

In order to derive the proof of Theorem \ref{l.bm}, i.e.,
the exact asymptotic behaviour of the infimum of the standard Brownian motion with drift
we shall investigate in \nelem{eq:diff} the tail asymptotics of the difference $X-Y$ assuming that $X$ has distribution $F$ with unbounded support and $Y\ge 0$ almost surely.
If for any $\eta>0$
\BQN\label{eta}
\limit{u} \frac{\pk{X> u + \eta}}{\pk{X>u}}=0,
\EQN
then  Lemma 2  in \cite{EmbHATHM} entails
\BQNY
\limit{u} \frac{ \pk{X- Y > u}}{\pk{X> u}} =   \pk{Y=0}.
\EQNY
If  $F$ is in the Gumbel max-domain of attraction 
with some positive scaling function $w(\cdot)$, i.e., 
\BQN\label{eqGU}
1- F(u+ x/ w(u)) = \exp(- x) (1- F(u)) (1+ o(1)), \quad \forall x\inr
\EQN
as $u\to\IF$, then \eqref{eta} is satisfied if additionally
 $\limit{u} w(u)= \IF$.
\COM{
 then, 
for any $\ve,\delta$ positive
\BQNY
 \pk{X- Y > u}&=& \pk{X-Y> u, Y\le \ve }+ \pk{X-Y> u, \ve + \delta  \ge Y> \ve }\\
 &&+ \pk{X-Y> u, Y> \ve+ \delta }\\
&\le & \pk{X-Y> u, Y\le \ve }+ \pk{X> u+ \ve}+ \pk{X> u+ \ve+ \delta }\\
&\le & \pk{X-Y> u, Y\le \ve }+ \pk{X> u+ \ve}(1+o(1)), \quad u\to \IF,
\EQNY
hence
\BQNY
 \pk{X- Y > u}&=&  \pk{X-Y> u, Y\le \ve }+ O(\pk{X> u}), \quad u\to \IF
\EQNY
and thus the probability that $Y$ takes small values is important. If $\pk{Y=0}>0$, then Lemma in \cite{EmbHATHM} implies
\BQNY
 \pk{X- Y > u}&=&  \pk{Y=0}\pk{X> u}(1+o(1)), \quad u\to \IF,
\EQNY
provided that $\limit{u} \pk{X> u + \eta}/\pk{X>u}=0$ for any $\eta>0$.\\
}
As shown below, it is possible to derive the exact tail asymptotics of $X-Y$ when $\pk{Y=0}=0$
assuming further  that for some $\alpha\ge 0$
\BQN\label{eqAL}
\pk{Y<  x/u}=x^\alpha \pk{Y < 1/u}\ooo, \quad \forall x> 0
\EQN
holds as $u\to \IF$.

\BEL\label{eq:diff}
 Let $X$ and $Y$ be two independent random variables. If \eqref{eqGU}
 holds for some positive function $w(\cdot)$ such that
$\limit{u} w(u)=\IF$ and further $Y\ge 0$ satisfies \eqref{eqAL} with some $\alpha\ge 0$, then we have
\BQN
\pk{X- Y > u} = \Gamma(\alpha+ 1) \pk{Y< 1/w(u )}\pk{X> u}\ooo, \quad u\to \IF.
\EQN
In particular, if $Y$ possesses a density function $f(\cdot)$ in a neighborhood of $0$ such that $f(0)>0$, then
\BQN
\pk{X- Y > u} = \frac{f(0)}{w(u)}\pk{X> u}\ooo, \quad u\to \IF.
\EQN
\EEL
\prooflem{eq:diff}
The assumption that $\limit{u} w(u)= \IF$ implies that $\exp(X)$ is in the Gumbel MDA with scaling function $w^*(u)= w(\ln u)/u$. Further, \eqref{eqAL} is equivalent with
$$ \limit{u} \frac{ \pk{e^{-Y}> 1- x/u}}{\pk{e^{-Y}> 1- 1/u}}= x^\alpha, \quad x>0. $$
Since for any positive $u$ we have
\BQNY
\pk{X- Y> u}&=& \pk{e^X e^{-Y}> e^u},
\EQNY
then by Example 1 in \cite{HashExt12} or Theorem 4.2 in \cite{EH15}
\BQNY
\pk{X- Y> u}&=& \Gamma(\alpha+1) \pk{ e^{-Y}> 1- 1/(e^u w^*(e^u))}\pk{e^{X} > e^{u}}\ooo \\
 &=& \Gamma(\alpha+1)  \pk{ Y< 1/ w(u)}\pk{e^{X}> e^{u}}\ooo, \quad u\to \IF.
\EQNY
In the special case that $Y$ possesses a density function $f(\cdot)$ with $f(0)>0$, then $\alpha=1$ and
$$\pk{ Y< 1/ w(u)} =  \frac{f(0)}{w(u)}\ooo$$
 as $u\to \IF$ establishing the proof. \QED

\COM{
\BEL \label{LemBm}
Let $\{B(t),t\ge0\}$ be a standard Brownian motion, and $c$ be a  nonnegative  constant. Then, for any positive constants $T_1,T_2$ satisfying $T_1<T_2$, we have
\BQN\label{eqLemBm}
\pk{\inf_{t\in[T_1,T_2]}\Bigl(B(t)-ct\Bigr)>u}=  \frac{K_{c,T_2-T_1}}{u} \pk{B(T_1)> u+c T_1}\ooo, \ u\to\IF,
\EQN
where
\BQN\label{eq:KcT}
K_{c,T_2-T_1}=2\varphi \LT( \frac{c(T_2-T_1)}{\sqrt{T_2-T_1}}\RT) \frac{1}{\sqrt{T_2-T_1}} -
 2 c  \Psi\LT( \frac{c(T_2-T_1)}{\sqrt{T_2-T_1}}\RT)>0.
 \EQN
\EEL

\prooflem{LemBm}
}
\tbb{{\it Proof of Theorem \ref{l.bm}}:}
Let $\NN$ be a standard $N(0,1)$ random variable \pE{with density function $\varphi$} which is independent of the Brownian motion
$B(\cdot):=B_1(\cdot)$. We have with $\Delta:= T_2- T_1>0$
\BQNY
\pk{\inf_{t\in[T_1,T_2]}\Bigl(B(t)-ct\Bigr)>u}
&=&\pk{\inf_{t\in[T_1,T_2]} \Bigl( B(t) - B(T_1) -c(t-T_1) + B(T_1) -cT_1\Bigr) >u}\\
&=&\pk{   T_1^{1/2} \NN  -\sup_{t\in[0,\Delta]}  \Bigl(B(t) +c t\Bigr) >u+c T_1}.
\EQNY
It is well-known that for any $u>0$  and $c\ge 0$ 
\BQNY
\pk{\sup_{t\in[0,\Delta]}  \Bigl(B(t) +c t\Bigr)> u} = \Psi\LT( \frac{u- c\Delta}{\sqrt{\Delta}}\RT)
+ e^{2 c   u} \Psi\LT( \frac{u+ c\Delta}{\sqrt{\Delta}}\RT),
\EQNY
\COM{consequently,
\BQNY
&& \pk{\inf_{t\in[T_1,T_2]}\Bigl(B(t)-ct\Bigr)>u}\\
&=& \int_{y\inr} \pk{  \sup_{t\in[0,\Delta]}  \Bigl(B(t) +c t\Bigr) <  u+ T_1^{1/2} y -c T_1  }\varphi(y)\, dy\\
&=& 1- \int_{y\inr} \pk{  \sup_{t\in[0,\Delta]}  \Bigl(B(t) +c t\Bigr) >  u+ T_1^{1/2} y -c T_1  }\varphi(y)\, dy\\
&=& 1- \int_{y\inr} \Bigl[ \Psi\LT( \frac{T_1^{1/2} y + u-c T_1  - c\Delta}{\sqrt{\Delta}}\RT)
+ e^{2 c   (T_1^{1/2} y -c T_1 -u )} \Psi\LT( \frac{T_1^{1/2} y -c T_1 -u + c\Delta}{\sqrt{\Delta}}\RT) \Bigr) \varphi(y)\, dy\\
&=& \int_{y\inr} \Phi\LT( \frac{T_1^{1/2} y +u -c T_2}{\sqrt{\Delta}}\RT) \varphi(y) \, dy
+  e^{-2 c   (c T_1 +u )} \int_{y\inr} e^{2 c   (T_1^{1/2} y -c T_1 -u )} \Psi\LT( \frac{T_1^{1/2} y -c T_1 -u + c\Delta}{\sqrt{\Delta}}\RT) \Bigr) \varphi(y)\, dy\\
\EQNY
}
hence the density function $q$  of $\sup_{t\in[0,\Delta]}  (B(t) +c t)$ is given by
\BQNY
q(u)=\varphi \LT( \frac{u- c\Delta}{\sqrt{\Delta}}\RT) \frac{1}{\sqrt{\Delta}}
- 2 c   e^{2 c  u} \Psi\LT( \frac{u+ c\Delta}{\sqrt{\Delta}}\RT)
+ e^{2 c   u} \varphi\LT( \frac{u+ c\Delta}{\sqrt{\Delta}}\RT)\frac{1}{\sqrt{\Delta}}, \quad u>0.
\EQNY
Since $\sqrt{T_1} \mathcal{N}$ has distribution in the Gumbel MDA with $w(u)=u/T_1$ and
$$ q(0)= 2\varphi \LT(  c \sqrt{\Delta} \RT) \frac{1}{\sqrt{\Delta}} -
 2 c   \Psi\LT(  c \sqrt{\Delta} \RT)>0 $$
the claim follows from \nelem{eq:diff}. \QED

\subsection{Proof of Theorem \ref{ThmGauBou}}\label{p.low.bound}
For any $u$ positive we have
\BQNY
\pk{\sup_{t\in[0,S]}\inf_{s\in[t,t+T]}\Bigl(X(s)-cs\Bigr)>u}
 &\ge&  \pk{\inf_{s\in[S,S+T]}\Bigl(X(s)-cs\Bigr)>u}\\
&=&\pk{\sup_{s\in[S,S+T]}\Bigl(-X(s)+cs\Bigr)<-u}.
\EQNY
Since we assume that $V(t):=\sigmaX^2(t)$ is a convex function and  $V(0)=0$, then for any $0\le s\le t$
\BQNY
V(t)\ge V(s)+V(t-s).
\EQNY
Therefore, by the Slepian lemma (e.g., \cite{Pit96})
\BQNY
\pk{\sup_{s\in[S,S+T]}\Bigl(-X(s)+cs\Bigr)<-u}&\ge&\pk{\sup_{s\in[S,S+T]}\Bigl(-B(V(s))+cs\Bigr)<-u}\\
&= &\pk{\inf_{t\in[S,S+T]} \Bigl(B(V(t))-ct\Bigr)>u}\\
&=&\pk{\inf_{t\in[V(S),V(S+T)]} \Bigl(B(t)-c g(t)\Bigr)>u},
\EQNY
where $B$ is a standard Brownian motion and $g(\cdot)$ is the inverse function of $V(\cdot)$. Further, since $g(s), s\ge0$ is differentiable, increasing and concave we have
(set $\rho_S=1/ V'(S)$ with $V'(t)$ the derivative of $V(t)$)
\BQNY
g(s)\le f(s):=\rho_S s+S-\rho_SV(S),\ \  s\ge0
\EQNY
implying thus 
\BQNY
\pk{\sup_{t\in[0,S]}\inf_{s\in[t,t+T]}\Bigl(X(s)-cs\Bigr)>u}
 &\ge&   \pk{\inf_{t\in[V(S),V(S+T)]} \Bigl(B(t)-c \rho_St\Bigr)>u+c(S-\rho_S V(S))}\\
 &=&K_{c\rho_S,V(S+T)-V(S)}\frac{V(S)}{u} \Psi\LT(\frac{ u+c S}{\sqrt{V(S)}}\RT)\ooo
\EQNY
as $u\to\IF$. where the last equality follows from \eqref{eqLemBm}, and $K_{c,y}$ is given as in \eqref{eq:KcT}.
\COM{
Next we show the asymptotical upper bound. We have, for any $b>0$
\BQNY
\pk{\sup_{t\in[0,S]}\inf_{s\in[t,t+T]}\Bigl(X(s)-cs\Bigr)>u}\le \pk{\sup_{t\in[0,S]}\inf_{s\in[t,t+bu^{-2/\alpha}]} \Bigl(X(s)-cs\Bigr)>u}.
\EQNY
Consequently, the claim follows directly from \netheo{ThmGauExact}.}
\QED

\subsection{Proof of Theorem \ref{th.log}}\label{p.log}

The proof follows straightforwardly from the combination of
Theorem \ref{ThmGauBou} and the fact that
\begin{eqnarray}
\mathcal{P}_S(u,T_u) \notag
&<&
\mathcal{P}_S(u)\\
&\le&
\pk{\sup_{t\in[0,S]} X(t)>u}\nonumber\\
&\le& \pk{\sup_{t\in[0,S]} B_1(\sigma^2(t))>u}\label{from.slepian}\\
&=&2\Psi\left(\frac{u}{\sigma(S)}\right),\nonumber
\end{eqnarray}
where (\ref{from.slepian}) follows from the Slepian  \pE{lemma} (recall that $\sigma^2(\cdot)$ is convex). \QED

\subsection{Proof of Theorem \ref{ThmGauExact}}\label{p.exact}

Let $\delta(u)=(\ln u/u)^{2/{\beta_1}}, u>0$ and set
\BQNY
\Pi(u)=\pk{\sup_{t\in[S-\delta(u),S]}\inf_{s\in[t,t+T_u]}\Bigl(X(s)-cs \Bigr)>u},\ \ u>0.
\EQNY
It follows that
\BQNY
\Pi(u)\le\Piter_S(u,T_u)=\pk{\sup_{t\in[0,S]}\inf_{s\in[t,t+T_u]}\Bigl(X(s)-cs \Bigr)>u}\le \Pi(u)+\Pi_o(u), 
\EQNY
where $\Pi_o(u)=\pk{\sup_{t\in[0,S-\delta(u)]}\Bigl(X(t)-ct \Bigr)>u}.$
We shall show that
\BQN\label{eq:Pio}
\Pi_o(u)=o(\Pi(u)),\ \ \ u\to\IF,
\EQN
which on the turn implies 
$$
\Piter_S(u,T_u)=\Pi(u)\ooo,\ \ \ u\to\IF.
$$
Next, we derive the exact tail asymptotics of $\Pi(u)$.  For notational simplicity we set
$$
g_u(t)=\frac{u+ct}{\sigmaX(t)},\  \ X_u(t)=\frac{X(t)}{\sigmaX(t)}\frac{g_u(S)}{g_u(t)},\ \ \cL{\sigma_{X_u}^2}(t)=\Var(X_u(t))\ \ t\ge0.
$$
By Assumption \textbf{A1} for any small $\vn\in(0,1)$, there exists some small $\theta>0 $ and $u_0>0$ such that
\BQN\label{eq:gu}
 \lefteqn{ (1-\vn)\frac{A}{\SigS}\abs{ t}^{\beta_1}I_{(t>0)}+(1\mp\vn)\frac{A_{\pm}}{\SigS}\abs{ t}^{\beta_2}I_{(t\le0)}}\nonumber\\
&\le& 1-\frac{g_u(S)}{g_u(S -t)}\\
&\le&(1+\vn)\frac{A}{\SigS}\abs{ t}^{\beta_1}I_{(t>0)}+(1\pm\vn)\frac{A_{\pm}}{\SigS}\abs{ t}^{\beta_2}I_{(t\le0)}\nonumber
\EQN
holds for all $t\in[-\theta,\theta]$ and all $u>u_0$. Note that in the derivation of the above inequality we used the fact that $\beta_1\le 1$ and ${\beta_2}\le1$.
By changing the time we obtain
\BQNY
\Pi(u)=\pk{\sup_{t\in[0,\delta(u)]}\inf_{s\in[0,T_u]}X_u(S+s-t) >g_u(S)}.
\EQNY
The idea for finding the exact asymptotics of $\Pi(u)$  is
analogous to the one used in \cite{Pit96}.
Let $q=q(u)=u^{-2/\alpha}$ and set for any $\lambda>T$
$$
\Del_k=\LT[k\lambda q,(k+1)\lambda q\RT],\  \cL{k\inn_0}, \ \ \text{and}\ \
 N(u)=\LT\lfloor \lambda^{-1}\delta(u)q^{-1}\RT\rfloor\cL{+1},
$$
 where $\lfloor\cdot\rfloor$ is the ceiling function.
 We shall  investigate separately the following three cases:

(i) $\alpha<{\beta_1}$,\ \ \ (ii) $\alpha={\beta_1}$,\ \ \ (iii) $\alpha>{\beta_1}$.

Since the case $T=0$ follows as a limiting result we shall consider for (i) and (ii) only $T\in(0,\IF)$.

\underline{(i) $\alpha<{\beta_1}$:}  We have by the Bonferroni inequality
\BQNY
 \sum_{k=0}^{N(u)}\pi_k(u)\ge\Pi(u)\ge  \sum_{k=0}^{N(u)-1}\pi_k(u)-\Sigma(u),
\EQNY
where
\BQNY
&&\pi_k(u)=\pk{\sup_{t\in \Del_k}\inf_{s\in[0,T_u]}X_u(S+s-t)>g_u(S)},\ \  \cL{k\inn_0},\\
&&\Sigma(u)=\underset{0\le i< j\le N(u)}{\sum\sum} \pk{\sup_{t\in \Del_i}\inf_{s\in[0,T_u]}X_u(S+s-t)>g_u(S), \sup_{t\in \Del_j}\inf_{s\in[0,T_u]}X_u(S+s-t)>g_u(S)}.
\EQNY
In view of \eqref{eq:gu}  for any $k=0,\cdots,N(u)$ 
\BQN\label{eq:sigxu}
\lefteqn{1-(1+\vn)\frac{A}{\SigS}\abs{ t-s}^{\beta_1}I_{(t>s)}-(1\pm\vn)\frac{A_{\pm}}{\SigS}\abs{ t-s}^{\beta_2}I_{(t\le s)}}\nonumber\\
&\le&\sigma_{X_u}(S+s-t)\\
&\le&1- (1-\vn)\frac{A}{\SigS}\abs{ t-s}^{\beta_1}I_{(t>s)}-(1\mp\vn)\frac{A_{\pm}}{\SigS}\abs{ t-s}^{\beta_2}I_{(t\le s)}\nonumber
\EQN
holds for all $(t,s)\in\Del_k\times[0,T_u]$. Define next
$$
Y_u(t,s)=\frac{X_u(S+s-t)}{\sigma_{X_u}(S+s-t)},\ \ t,s\in[0,S].
$$
For any small $\vn\in(0,1)$ and $k=1,\cdots, N(u)$
\BQNY
\pi_k(u)\le\pk{\sup_{t\in \Del_k}\inf_{s\in[0,T_u]}Y_u(t,s)>g_u(S)\LT(1+(1-\vn)^2\frac{A}{\SigS}\abs{ k\lambda q-T_u }^{\beta_1}\RT)}
\EQNY
and
\BQNY
\pi_k(u)\ge\pk{\sup_{t\in \Del_k}\inf_{s\in[0,T_u]}Y_u(t,s)>g_u(S)\LT(1+(1+\vn)^2\frac{A}{\SigS}\abs{ (k+1)\lambda q}^{\beta_1}\RT)}
\EQNY
are valid for $u$ sufficiently large. 
Moreover, for $u$ sufficiently large also
\BQNY
\pi_0(u)\le\pk{\sup_{t\in \Del_0}\inf_{s\in[0,T_u]}Y_u(t,s)>g_u(S)\LT(1+(1\mp\vn)^2\frac{A_\pm}{\SigS}\abs{f_\pm(u) }^{\beta_2}\RT)}
\EQNY
and
\BQNY
\pi_0(u)\ge\pk{\sup_{t\in \Del_0}\inf_{s\in[0,T_u]}Y_u(t,s)>g_u(S)\LT(1+(1+\vn)^2\frac{A}{\SigS}\abs{\lambda q}^{\beta_1}+ (1\pm\vn)^2\frac{\abs{A_\pm}}{\SigS}\abs{h_\pm(u) }^{\beta_2}\RT)}
\EQNY
are valid, where $f_+(u)=h_-(u)=0$, $f_-(u)=T_u+\lambda q$ and $h_+(u)=T_u$.
Consequently, an application of \nelem{LemGP} in Appendix yields that
\BQNY
&& \sum_{k=1}^{N(u)}\pi_k(u)\le   \sum_{k=1}^{N(u)}\pk{\sup_{t\in [0,\lambda]}\inf_{s\in[0,T]}Y_u(t q+k\lambda q,sq)>g_u(S)\LT(1+(1-\vn)^2\frac{A}{\SigS}\abs{k\lambda q-T_u}^{\beta_1}\RT)}\\
 &&=\FF_\alpha(\hat a \lambda, \hat a T)\frac{1}{\sqrt{2\pi} g_u(S)} \sum_{k=1}^{N(u)}\exp\LT(-\frac{(g_u(S))^2\LT(1+(1-\vn)^2\frac{A}{\SigS}\abs{k\lambda q-T_u}^{\beta_1}\RT)^2}{2}\RT)\ooo
\EQNY
as $u\to \IF$, where $\hat a=D^{1/\alpha} \SigS^{-2/\alpha}$. Further, since
$$
\int_0^\IF\exp(-bx^{\beta_1})dx=\Gamma\LT(\frac{1}{{\beta_1}}+1\RT)b^{-\frac{1}{{\beta_1}}},\ \ b>0,{\beta_1}>0
$$
we have
\BQNY
 \sum_{k=1}^{N(u)}\pi_k(u)\le   \frac{1}{\lambda}\FF_\alpha(\hat a \lambda, \hat a T) \Gamma\LT(\frac{1}{{\beta_1}}+1\RT)\LT(\frac{\SigS^3}{(1-\vn)^2A}\RT)^{\frac{1}{{\beta_1}}} u^{\frac{2}{\alpha}-\frac{2}{{\beta_1}}}\Psi(g_u(S)) \ooo
\EQNY
as $u\to\IF.$ Similarly
\BQNY
 \sum_{k=1}^{N(u)-1}\pi_k(u)\ge   \frac{1}{\lambda}\FF_\alpha(\hat a \lambda, \hat a T) \Gamma\LT(\frac{1}{{\beta_1}}+1\RT)\LT(\frac{\SigS^3}{(1+\vn)^2A}\RT)^{\frac{1}{{\beta_1}}} u^{\frac{2}{\alpha}-\frac{2}{{\beta_1}}}\Psi(g_u(S)) \ooo
\EQNY
as $u\to\IF.$ By  \nelem{LemGP} and our assumption $\alpha<\beta_1\le\beta_2$ we obtain
$$
\pi_0(u)=\FF_\alpha(\hat a \lambda, \hat a T)\Psi(g_u(S))\ooo=o\LT(\sum_{k=1}^{N(u)-1}\pi_k(u)\RT)
$$
as $u\to\IF.$ Further, we have (set $\theta_i(u):= 1+(1-\vn)^2\frac{A}{\SigS}\abs{\max(0,i\lambda q-T_u)}^{\beta_1}+(1\mp\vn)^2\frac{A_\pm}{\SigS}\abs{f_\pm(u) }^{\beta_2}$)
\BQNY
\Sigma(u)\le\underset{0\le i< j\le N(u)}{\sum\sum} \pk{\sup_{t\in \Del_i} Y_u(t,0)>g_u(S)\theta_i(u),
\sup_{t\in \Del_j}Y_u(t,0)>g_u(S) \theta_i(u)}.
\EQNY
Letting $\vn\to0$ and $\lambda\to\IF$ we conclude by similar arguments as in the proof of Theorem 3.1 in \cite{DHJ13b} that
\BQNY
\Pi(u) =\FF_\alpha(\hat a T)  \Gamma\LT(\frac{1}{{\beta_1}}+1\RT)D^{\frac{1}{\alpha}} A^{-\frac{1}{{\beta_1}}}   \SigS^{\frac{3}{{\beta_1}}-\frac{2}{\alpha}} u^{\frac{2}{\alpha}-\frac{2}{{\beta_1}}}\Psi(g_u(S)) \ooo
\EQNY
as $u\to \IF$, and $\FF_\alpha(T)\in(0,\IF)$. 

\underline{(ii) $\alpha={\beta_1}$:} We use the same notation as in Case (i). We have by the Bonferroni inequality
\BQNY
 \pi_0(u)\le\Pi(u)\le \pi_0(u)+ \sum_{k=1}^{N(u)}\pi_k(u).
\EQNY
It follows from \nelem{LemGP} that
\BQN\label{eq:pi0}
\pi_0(u)=\PP_{\alpha,\beta_2}^{A/(D\SigS), A_\pm/(D\SigS)}(\hat a \lambda, \hat a T)\Psi(g_u(S))\ooo, \ u\to\IF.
\EQN
Further, for any small $\vn\in(0,1)$
\BQNY
\sum_{k=1}^{N(u)}\pi_k(u)\le\sum_{k=1}^{N(u)}\pk{\sup_{t\in \Del_0}Y_u(t+k\lambda q,0)>g_u(S)\LT(1+(1-\vn)^2\frac{A}{\SigS}\abs{k\lambda q-T_u}^{\beta_1}\RT)}
\EQNY
for  $u$ sufficiently large. Using \nelem{LemGP} (or Lemma 1 in \cite{DebKo2013}) we have further that
\BQNY
\sum_{k=1}^{N(u)}\pi_k(u)\le G \FF_\alpha(\hat a \lambda,0) \Psi(g_u(S)) \sum_{k=1}^{\IF}\exp\LT(-\frac{A}{2 \SigS^3}(k\lambda-T)^{\beta_1}\RT)\ooo
\EQNY
as $u\to\IF$, for some positive constant $G$.
Therefore, we conclude that, for any $\lambda_1,\lambda_2>T$
\BQNY
&&\PP_{\alpha,\beta_2}^{A/(D\SigS), A_\pm/(D\SigS)}(\hat a \lambda_2, \hat a T)\le\liminf_{u\to\IF}\frac{\Pi(u)}{\Psi(g_u(S))}\le\limsup_{u\to\IF}\frac{\Pi(u)}{\Psi(g_u(S))}\\
&&\le \PP_{\alpha,\beta_2}^{A/(D\SigS), A_\pm/(D\SigS)}(\hat a \lambda_1, \hat a T)+G \FF_\alpha(\hat a \lambda_1,0) \sum_{k=1}^{\IF}\exp\LT(-\frac{A}{2 \SigS^3}(k\lambda_1-T)^{\beta_1}\RT).
\EQNY
Further, it follows from Corollary D.1 in \cite{Pit96} that $\FF_\alpha(\hat a \lambda_1,0)=\H_\alpha(\hat a \lambda_1)\le \lfloor \hat a\lambda_1\rfloor+1$, and thus
$$
\lim_{\lambda_1\to\IF}\FF_\alpha(\hat a \lambda_1,0) \sum_{k=1}^{\IF}\exp\LT(-\frac{A}{2 \SigS^3}(k\lambda_1-T)^{\beta_1}\RT)=0.
$$
Consequently, by letting $\lambda_1$ and $\lambda_2$ tend to infinity, respectively, we conclude that
\BQNY
 \lim_{u\to\IF}\frac{\Pi(u)}{\Psi(g_u(S))}=\PP_{\alpha,\beta_2}^{A/(D\SigS), A_\pm/(D\SigS)}( \hat a T)\in(0,\IF).
\EQNY
\underline{(iii) $\alpha>{\beta_1}$:} We use the same notation as in Case (i) and Case (ii). 
In view of \eqref{eq:gu} and the fact that $\lim_{u\to\IF}T_uu^{2/\alpha}= 0$, for any small $\vn, \vn_1\in(0,1)$
\BQNY
\Pi(u)& \ge &
\pk{ \inf_{s\in[0,T_u]}X_u(S+s) >g_u(S)}\\
&\ge& \pk{ \inf_{s\in[0,T_u]}Y_u(0,s) > g_u(S) \LT(1+(1+\vn)\frac{\abs{A_\pm}}{\SigS}T_u^{\beta_2}\RT)}\\
&\ge&\pk{ \inf_{s\in[0,\vn_1]}Y_u(0,su^{-\frac{2}{\alpha}}) > g_u(S) \LT(1+(1+\vn)\frac{\abs{A_\pm}}{\SigS}T_u^{\beta_2}\RT)}
\EQNY
holds for all $u$ sufficiently large. Moreover, it follows from \nelem{LemGP}  that
\BQNY
&&\pk{ \inf_{s\in[0,\vn_1]}Y_u(0,su^{-\frac{2}{\alpha}}) > g_u(S) \LT(1+(1+\vn)\frac{\abs{A_\pm}}{\SigS}T_u^{\beta_2}\RT)}\\
&&= \H_\alpha^{\inf}(\hat a\vn_1) \Psi\LT(g_u(S) \LT(1+(1+\vn)\frac{\abs{A_\pm}}{\SigS}T_u^{\beta_2}\RT)\RT)\ooo
\EQNY
as $u\to \IF,$ where
$$
\H_\alpha^{\inf}(T)=\E{\exp\LT(\inf_{t\in[0,T]}\LT(\sqrt{2}B_\alpha( t )-t^{\alpha}\RT) \RT)},\ \ T\ge0.
$$
  Therefore, letting $\vn,\vn_1\to 0$ we have by the fact that $\lim_{u\to\IF}T_u  u^{2/{\beta_2}}= 0$
\BQNY
\Pi(u)&\ge& \Psi\LT(g_u(S) \LT(1+\frac{\abs{A_\pm}}{\SigS}T_u^{\beta_2}\RT)\RT)\ooo\\
 &=& \Psi (g_u(S))\ooo 
\EQNY
as $u\to \IF.$
Next we give the upper bound.  Since $\alpha>{\beta_1}$, we have
\BQNY
\Pi(u)\le \pk{\sup_{t\in \Del_0} X_u( S+T_u-t)>g_u(S)}.
\EQNY
Further, 
for any small $\vn\in(0,1)$,
 \BQNY
\Pi(u)&\le&\pk{\sup_{t\in \Del_0} Y_u(t,T_u)>g_u(S)\LT(1-(1+\vn)\frac{\abs{A_\pm}}{\SigS}T_u^{\beta_2}\RT) }\\
&=& \H_\alpha^{\sup}(\hat a\lambda) \Psi\LT(g_u(S)\LT(1-(1+\vn)\frac{\abs{A_\pm}}{\SigS}T_u^{\beta_2}\RT)\RT)\ooo
\EQNY
as $u\to \IF,$ where the last equation follows from \nelem{LemGP}, and
$$
\H_\alpha^{\sup}(T)=\E{\exp\LT(\sup_{t\in[0,T]}\LT(\sqrt{2}B_\alpha( t )-t^{\alpha}\RT) \RT)},\ \ T\ge0.
$$
Consequently, letting $\lambda,\vn\to 0$
\tbb{and using that $\lim_{u\to\IF}T_u  u^{2/{\beta_2}}= 0$,
we conclude
that, as $u\to \IF$,}
\BQNY
\Pi(u)  \le  \Psi (g_u(S))\ooo. 
\EQNY
\tbb{Thus} the claim follows.

\underline{Proof of \eqref{eq:Pio}.} First, for any fixed $\vn\in(0,1)$, by {\bf A1} we can choose some small $\theta_0>0$ such that
$$
\sigmaX(t)\le \SigS-(1-\vn)A(S-t)^{\beta_1}
$$
holds for all $t\in[S-\theta_0, S]$. Additionally, this $\theta_0$ can also be chosen such that
$$
\sup_{t\in [0,S-\theta_0)}\sigmaX(t)<\sigmaX(S-\theta_0)<\SigS.
$$
Clearly,
\BQNY
\Pi_o(u)\le\pk{\sup_{t\in[0,S-\theta_0]} \Bigl(X(t)-ct \Bigr)>u}+\pk{\sup_{t\in[S-\theta_0,S-\delta(u)]}\Bigl( X(t)-ct \Bigr)>u}=:\Pi_1(u)+\Pi_2(u).
\EQNY
By Borell-TIS inequality (cf. \cite{Pit96})
\BQNY
\Pi_1(u)\le\pk{\sup_{t\in[0,S-\theta_0]} X(t)>u}\le \exp\LT(-\frac{\LT(u-\E{\sup_{t\in[0,S]}X(t)}\RT)^2}{2 \sigmaX^2(S-\theta_0)}\RT)
\EQNY
for $u$ sufficiently large. Further, by \textbf{A3} we have applying Theorem 8.1 in \cite{Pit96}
\BQNY
\Pi_2(u)&\le&\pk{\sup_{t\in[\delta(u),\theta_0]} X(S-t)>u}\\
&\le& G u^{\frac{2}{\gamma}+1}\exp\LT(-\frac{u^2}{2 \SigS^2}\LT(1+(1-\vn)\frac{A}{\SigS}(\delta(u))^{\beta_1}\RT)\RT)
\EQNY
for $u$ sufficiently large, where $G$ is some positive constant independent of $u$.  Consequently,
we conclude from the asymptotics of $\Pi(u)$ for all the cases above that $\Pi_o(u)= o(\Pi(u))$, and thus the proof is complete. \QED 

\subsection{Proof of Corollary \ref{CorfBm}}
Since $X$ is a fBm with Hurst index $\alpha/2$ we have that
\BQNY
\sigmaX(t)=t^{\frac{\alpha}{2}}=S^{\frac{\alpha}{2}}-\frac{\alpha}{2} S^{{\frac{\alpha}{2}}-1}(S-t)\ooo,\ \ \ t\to S,
\EQNY
and
\BQNY
 \Cov\LT(\frac{X(t)}{\sigmaX(t)}, \frac{X(s)}{\sigmaX(s)}\RT)=1- \frac{1}{2S^\alpha}|t-s|^{\alpha}\ooo, \quad  t,s \to S.
 \EQNY
Moreover, for any $s,t\ge0$
\BQNY
\E{(X(t)-X(s))^{2}}  =   |t-s|^{\alpha}.
\EQNY
Consequently, the claim follows by an application of \netheo{ThmGauExact}. \QED

\COM{
\prooftheo{ThmRuintime} Denote $S_u=S-xu^{-2}$. We have,  for any $x\ge0$
\BQNY
\pk{u^2(S-\tau_u)\le x\big\lvert \tau_u<S}=\frac{\pk{\tau_u<S_u}}{\pk{\tau_u<S}}=\frac{\Piter_{S_u}(u,T_u)}{\Piter_{S_u}(u,T_u)}. 
\EQNY
As in the proof of \netheo{ThmGauExact} let $g_u(t)=\frac{u+ct}{\sigmaX(t)}$, and further define $Z_u(t)=\frac{X(t)}{\sigmaX(t)}\frac{g_u(S_u)}{g_u(t)}$. Then
\BQNY
\Piter_{S_u}(u,T_u)=\pk{\sup_{t\in[0,1]}\inf_{s\in[t,t+T_uS_u^{-1}]} Z_u(S_u s)>g_u(S_u)}.
\EQNY
}

\section{Appendix}
\def\bD{ \mathbf{ D} }
\def\sigxiu{\sigma_{\xi_u}}

Let $\bD $ be a compact set in $\R^n$, $n\inn$ and suppose without loss of generality that $\vk{0}\in \bD$.
Further, let $\{\xi_u(\vk{t}), \vk{t}\in\bD\}$, $u>0$
be a family of centered Gaussian random fields with a.s. continuous sample paths and variance function $\sigma_{\xi_u}^2(\vk{\cdot})$. Below $||\cdot||$ stands for the Euclidean norm in $\R^n$.  We assume that $\xi_u$ satisfies the following conditions:

{\bf C1:} 
 $\sigxiu(\vk{0})=1$ for all  $u$ large, and there exists some  \cL{bounded measurable} function $d(\vk{\cdot})$ on $\bD$  such that
$$\lim_{u\rw\IF}\sup_{\vk{t}\in\bD}
\abs{u^2(1-\sigxiu(\vk{t}))-d(\vk{t})}=0.$$

{\bf C2:} There exist some centered Gaussian random field $\{\eta(\vk{t}), \vk{t}\in\R^n\}$ with a.s. continuous sample paths, $\eta(\vk{0})=0$ and variance function $\sigma_\eta^2(\vk{\cdot})$ such that
$$
\lim_{u\rw\IF}
 u^2\Var(\xi_u(\vk{t})-\xi_u(\vk{s}))=2\Var(\eta(\vk{t})-\eta(\vk{s}))
 $$
  holds for all $\vk{t},\vk{s}\in  \bD$.

 {\bf C3:} There exist some  constants $G, \nu>0, u_0>0,$ such that, for any $u>u_0$  
\BQNY
u^2\Var(\xi_u(\vk{t})-\xi_u(\vk{s})) \le G \abs{ \abs{\vk{t}-\vk{s}}}^{\nu}
\EQNY
holds uniformly with respect to $\vk{t},\vk{s}\in  \bD$.

As in \cite{DebKo2013} let $F: C(\bD)\to \R$   be a continuous functional acting on $C(\bD)$, the space of continuous functions on the compact set $\bD$. Assume that:

{\bf F1:} $\abs{F(f)}\le \sup_{\vk{t}\in\bD} \pE{\abs{f(\vk{t})}}$ for any $f\in C(\bD)$.

{\bf F2:} $F(af+b)=aF(f)+b$ for any $f\in C(\bD)$ and $a>0,b\inr$.

For any \cL{bounded measurable} function $d(\vk{\cdot})$ on $\bD$ with $d(\vk{0})=0$ and  $F$ satisfying {\bf F1} we define a constant
\BQN\label{eq:GPitCons}
\H_{\eta,d}^F(\bD)=\E{\exp\LT(F\LT(\sqrt{2}\eta(\vk{t})-\sigma_\eta^2(\vk{t})-d(\vk{t})\RT)  \RT)}.
\EQN
Along the lines of the proof in \cite{DebKo2013} we get that $\H_{\eta,d}^F(\bD)\in(0,\IF)$.

\tbb{The following result generalizes} Lemma 6.1 in \cite{Pit96} and Lemma 1 in \cite{DebKo2013}.
\BEL\label{LemGP}
Let $\{\xi_u(\vk{t}), \vk{t}\in\bD\}$, $u>0$
be the family of centered Gaussian random fields defined as above satisfying {\bf C1-C3} with some function
 $d(\vk{\cdot})$ and some Gaussian random field $\eta$. Let $F: C(\bD)\to \R$   be a continuous functional such that  {\bf F1-F2} hold.
 Then, for any positive  measurable  function $g(\cdot)$   satisfying $\lim_{u\rw\IF}g(u)/u=a \in (0,\IF)$
\BQN\label{eq:lem1}
\pk{F(\xi_u)>g(u)}=\H_{a\eta,a^2d}^F(\bD)\Psi(g(u))\oo
\EQN
holds as $u\to\IF$, provided that $\pk{F(\xi_u)>g(u)}>0$ for all large $u$.
\EEL
\def\xiu{\xi_u}

\prooflem{LemGP} The proof is based on the classical approach rooted in the ideas of \cite{PicandsA, Pit96}. For all
$u>0$ large
\BQN\label{eq:lamS1}
\ \pk{F(\xi_u)>g(u)}=\Psi(g(u))\int_{\R} \exp\LT(w-\frac{w^2}{2(g(u))^2}\RT)\pk{F(\xi_u)>g(u)\Bigl| \xi_u(\vk{0})=g(u)-\frac{w}{g(u)}}dw.
\EQN
Let, for any $u>0, w\inr$,
$\zeta_u=\{\zeta_u(\vk{t})=g(u)(\xi_u(\vk{t})-g(u))+w,   \vk{t}\in\bD\}.$
Using {\bf F2} the conditional probability in the integrand of \eqref{eq:lamS1} can be written as
\BQNY
 \pk{F(\xi_u)>g(u)\Bigl| \xi_u(\vk{0})=g(u)-\frac{w}{g(u)}}=\pk{F(\chi_u)>w},
\EQNY
where $\chi_u=\zeta_u|\zeta_u(\vk{0})=0$.
Denote
\BQNY
R_{\xiu}(\vk{t},\vk{s})=\E{\xiu(\vk{t})\xiu(\vk{s})},\ \ \vk{s},\vk{t} \in\bD
\EQNY
to be the covariance function of $\xi_u$. We have that
the conditional random field
$\chi_u=
\LT\{\chi_u(\vk{t}) , \vk{t}\in\bD\RT\}
$
has the same finite-dimensional distributions as
\BQNY
\LT\{g(u)(\xiu(\vk{t})-R_{\xiu}(\vk{t},\vk{0})\xiu(\vk{0}))-(g(u))^2(1-R_{\xiu}(\vk{t},\vk{0}))+w(1-R_{\xiu}(\vk{t},\vk{0})), \vk{t}\in\bD\RT\}.
\EQNY

Therefore, the following convergence
\BQNY
\E{\chi_u(\vk{t})}=-(g(u))^2(1-R_{\xiu}(\vk{t},\vk{0}))+w(1-R_{\xiu}(\vk{t},\vk{0})) \to -a^2(\sigma_\eta^2(\vk{t})+d(\vk{t})), \  u\rw\IF 
\EQNY
holds, for any $w\in\R$, uniformly with respect to $\vk{t}\in\bD.$ Moreover, for any $\vk{t},\vk{s}\in\bD$ we have
\BQNY
\Var\Bigl( \chi_u(\vk{t})-\chi_u(\vk{s})\Bigr)&=&(g(u))^2\LT(\E{\Bigl( \xiu(\vk{t})-\xiu(\vk{s})\Bigr)^2}-\LT(R_{\xiu}(\vk{t},\vk{0})-R_{\xiu}(\vk{s},\vk{0})\RT)^2\RT)\\
& \to & 2a^2\Var(\eta(\vk{t})-\eta(\vk{s})), \  u\rw\IF.
\EQNY
Therefore, the finite-dimensional distributions of $\chi_u$ converge to those of $\widetilde{\eta}=\{\sqrt{2}a \eta(\vk{t})-\sigma_{a\eta}^2(\vk{t})-a^2d(\vk{t}), \vk{t}\in\bD\}$, whereas the tightness follows by Proposition 9.7 in \cite{Pit20}.
 The rest of the proof repeats  line-by-line
that of Lemma 1 in \cite{DebKo2013}.
\QED

\COM{
{\bf Example.} $X$ is $N(0,1)$ then $w(u)=u$. If $Y$ is regularly varying at zero with index 1 (for instance $Y= \abs{Z}$ where $Z$ is an $N(0,1)$ rv, then we have
\BQNY
\pk{X- Y > u} \sim \Gamma(2) \pk{X> u} \pk{Y< 1/u }, \quad u\to \IF.
\EQNY
So if $Y= \abs{Z}$ we get
\BQN
\pk{X- Y > u} \sim \Gamma(2) \pk{X> u} \frac{1}{ \sqrt{ 2 \pi} u} \sim \frac{1}{2 \pi u^2} \exp(- u^2/2) , \quad u\to \IF.
\EQN
}
\bigskip

{\bf Acknowledgement}: We are grateful to the editor and the referees  for their helpful suggestions. We also kindly acknowledge partial support by
the Swiss National Science Foundation Grants 200021-140633/1
\tbb{and the project RARE -318984 (an FP7 Marie Curie IRSES Fellowship)}.
KD also acknowledges partial support by NCN Grant No 2013/09/B/ST1/01778 (2014-2016).

\bibliographystyle{plain}

 \bibliography{gR}
\end{document}